\documentclass[11pt,reqno]{article}
\usepackage{amsmath}
\usepackage{mathrsfs}
\usepackage{amsfonts}
\usepackage{amssymb}

\usepackage{amssymb}
\usepackage{amsthm}
\usepackage{graphicx}              
\usepackage{amsmath}               
\usepackage{amsfonts}              
\usepackage{amsthm}                
\usepackage{setspace}
\usepackage{epstopdf}
\usepackage{epsfig}
\usepackage[]{xcolor}

\newtheorem{Theorem}{Theorem}[section]
\newtheorem{Proposition}{Proposition}[section]
\newtheorem{Lemma}{Lemma}[section]
\newtheorem{Corollary}{Corollary}[section]
\newtheorem{remark}{Remark}[section]

\renewcommand{\thefootnote}{\fnsymbol{footnote}}
\makeatletter

\newcommand{\Rmnum}[1]{\expandafter\@slowromancap\romannumeral #1@}
\makeatother

\allowdisplaybreaks
\numberwithin{equation}{section}

\title{
Decay and scattering of solutions to nonlinear Schr\"odinger equations with regular potentials for nonlinearities of sharp growth}

\author{ Ze Li \qquad Lifeng Zhao }

\date{}

\begin{document}

\maketitle

\renewcommand{\thefootnote}{\fnsymbol{footnote}}

\noindent {\bf{Abstract } }
In this paper, we prove the decay and scattering in the energy space for nonlinear Schr\"odinger equations with regular potentials in $\Bbb R^d$ namely, $i{\partial _t}u + \Delta u - V(x)u + \lambda |u|^{p - 1}u = 0$. We will prove decay estimate and scattering of the solution in the small data case when $1+\frac{2}{d}<p\le1+\frac{4}{d-2}$, $d\ge3$. The index $1+\frac{2}{d}$ is sharp for scattering concerning the result of W. Strauss \cite{W2}.

\noindent{\bf MR Subject Classification:} 35Q55.

\section{Introduction}\label{se1}

In this paper, we consider the nonlinear Schr\"odinger equation with a potential:
\begin{align}\label{1}
\left\{ \begin{array}{l}
 i{\partial _t}u + \Delta_V u + \lambda |u{|^{p - 1}}u = 0, \\
 u(h,x) = {u_0}(x). \\
 \end{array} \right.
\end{align}
where $u:[h,\infty)\times{ \Bbb R}^d\to {\Bbb C}, h> 0$, $\Delta_V=\Delta-V, V: \mathbb{R}^d \to \mathbb{R}$, $\lambda= \pm 1$ and $1< p < \infty$.
When $ p = 1 + \frac4d,\, d\ge 1$ and $1+ \frac4{d-2},\, d\ge 3$, the equation is called mass-critical and energy-critical respectively. The equation is called mass-supercritical for $d\ge 1$ if $p >  1 + \frac4d $
and energy-subcritical for $d\ge 3$ if $p < 1+\frac4{d-2}$.  If $\lambda=-1$, the equation is called defocusing; otherwise, it is called focusing if $\lambda = 1$.

There are many important areas of application which motivate the study of nonlinear Schr\"odinger equations with potentials (Gross-Pitaevskii equation). In the most fundamental level, it arises as a mean field limit model governing the interaction of a plenty large number of weakly interacting bosons \cite{Hepp,Lieb,Spohn}. In a macroscopic level, it arises as the equation governing the evolution of the envelope of the electric field of a light pulse propagating in a medium with defects, see for instance, \cite{GSW,GWH}.

In this paper, we aim to prove scattering results for (\ref{1}) as what has been done in the case of nonlinear Schr\"odinger equations without potentials. More precisely, we will show the scattering in the small data case for any $1+\frac{2}{d}<p\le1+\frac{4}{d-2}$, $d\ge3$, which is sharp concerning the result of Strauss \cite{W2}.

The decay and scattering for small initial data has been studied for decades. When $V=0$, it has been shown that for $d\ge 1,\ p=1+\frac{2}{d}$ is the critical exponent for scattering. In fact, for $1 + \frac{2}{d} < p < 1 + \frac{4}{d},\ d\ge 1$, decay and scattering of the solution in the small data case was proved by McKean and Shatah \cite{MS}. For $1 + \frac{4}{d} \le p \le 1+ \frac4{d-2}, \ d\ge 3$ and $1 + \frac{4}{d} \le p < \infty,\ d= 1,2$, local wellposedness and small data scattering was proved by Strauss \cite{S2}, we also refer to \cite{TC}. Moreover,  Strauss \cite{W2} showed when $1< p\le 1+\frac2d$ for $d\ge 2$ and $1< p \le 2$ for $d=1$, the only scattering solution is zero. This was extended to the case $1< p \le 3$ for $d=1$ by Barab \cite{B}. The existence and the form of the scattering operator was obtained by Ozawa \cite{O} for $d=1$ and by Ginibre and Ozawa \cite{GO} for $d\ge2$. The completeness of the scattering operator and the decay estimate were obtained by Hayashi and Naumkin \cite{HN}. For all solutions, not only for small ones, for $d=1$ in defocusing case, the completeness of the scattering operator and decay were obtained by Deift and Zhou \cite{DZ}.

When $V\ne 0$, the situation is much more involved. In \cite{SVN}, Cuccagna, Georgiev and Visciglia proved decay and scattering for small initial data for $p>3$ in one dimension when the potential $V$ is real Schwartz function with $\sigma(-\Delta_V) = [0,\infty)$.

In the article, we will consider the potentials satisfying the following assumptions:

\noindent{\bf{ Regular Potential Hypothesis}}\\
Suppose that $V$ is a real-valued potential satisfying\\
{\it (i)} ${\left\langle x \right\rangle ^{2N+1}}\left( {\left| V(x) \right| + \left| {\nabla V}(x) \right|} \right) \in {L^\infty(\mathbb{R}^d) }$, for some $ N>d$;\\
{\it (ii)} the spectrum of $-\Delta_V$ is continuous, and 0 is neither a resonance nor an eigenvalue of $-\Delta_V$;\\
{\it (iii)} ${\left\langle x \right\rangle ^\alpha }V(x)$ is a bounded operator from ${H^\eta }$ to ${H^\eta }$ for some $\alpha  > d + 4$, $\eta > 0$ with $\mathcal{F}{V}\in L^1$;

Define  $\|f\|_{\Sigma}=\|\left\langle {x} \right\rangle^2  f\|_2+\|\left\langle {x} \right\rangle \nabla f\|_2+\|f\|_{H^2}.$ The first main theorem is as follows which extends the results of \cite{SVN} to higher dimensions($d\ge3$).

\begin{Theorem}\label{th5}
Suppose that $V$ satisfies the Regular Potential Hypothesis and its Kato norm satisfies
\begin{equation}\label{condition}
\sup \limits_x\int \frac{|V(y)|}{|x-y|^{d-2}}\,\mathrm{d}y  <\frac{1}{c(d)},
\end{equation}
where $c(d)$ is a constant in the estimate of the Bessel function of the third type and $c(3)=\frac{1}{2\pi}$.
For $d\ge3$, $1+\frac{2}{d}<p \le 1+\frac{4}{d-2}$, $\lambda=\pm1$, if the initial data $\|u_0\|_{\Sigma}$ is sufficiently small, then \eqref{1} is globally well-posed. Moreover, for any $\gamma< \frac{d}2\big( 1  - \frac{1}{{p}}\big)$, we have the following decay estimate:
\begin{equation}\label{eq1.3}
\|u(t,x)\|_{L^{2p}_x}\le Ct^{ -\gamma},
\end{equation}
and as a consequence of the decay estimate, there exists $u_+\in H^1$ such that
\begin{align*}\label{ko}
\mathop {\lim }\limits_{t\to \infty}\|u(t)-e^{it\Delta_V}u_+\|_{H^1}=0.
\end{align*}
\end{Theorem}

\begin{remark}\label{re1.2}
The assumption {\it (iii)} is used to get a dispersive estimate of $e^{it\Delta_V}$. The condition given here is due to  Journe, Soffer and Sogge \cite{JAC}. There are many related works in this direction such as \cite{RS,BG}.
\end{remark}

The proof of Theorem \ref{th5} is based on the commutator method introduced in \cite{SVN}.
In fact, scattering can be reduced to decay estimates, namely
\begin{equation}\label{2}
\|u(t,x)\|_{L_x^{\infty}}\le Ct^{-\frac{d}2}.
\end{equation}
Then by introducing a vector field $|J_V|^su$ which is roughly $t^s{(-\Delta)}^{\frac{s}2}u$, it suffices to prove
\begin{equation}\label{31}
\|u(t,x)\|_{L_x^\infty}\le Ct^{-s}{\big\||J_V|^su\big\|^{\theta}_{L^{\infty}_tL^2_x}}\|u\|^{1-\theta}_2.
\end{equation}
Moreover, $|J_V|^su$ satisfies
$$i\partial_t|J_V|^su+\Delta_V |J_V|^su + B(s)+|J_V|^s\big(|u|^{p-1}u\big)=0.$$
From the Sobolev embedding, in order to realize (\ref{31}), $s$ should be chosen larger than $\frac{d}2$.
However, when $s\ge2$, the $B(s)$ term in the above equation is too complicated to handle. In order to overcome this difficulty, we introduce some new ideas. First, we find that the $L^{2p}$ decay estimate, roughly
\begin{equation}\label{3}
\|u\|_{2p}\le Ct^{-\frac{d}2 (1-\frac{1}{p})}\big\||J_V|^su\big\|_{L^{\infty}_tL^2_x},
\end{equation}
where $s= \frac{d}2( 1 -\frac{1}{ p})$, is enough for scattering. Since $s<2$ for $p\le1+\frac4{d-2}$, each term in the equation of $|J_V|^su$ can be estimated easily in our case.
Second, to establish (\ref{3}), we translate it to the corresponding estimate of the inverse of $|J_V|^s$, which can be reduced to the $L^p$ estimate of resolvent.
Third, we observe the ``almost equivalence" lemma, which is Lemma \ref{34} in the article, is sufficient to prove global well-posedness and scattering.

The article is organized as follows. In Section \ref{se2}, we prove resolvent estimates for the Schr\"odinger operators. Section \ref{se3} is devoted to the proof of Theorem \ref{th5}.
\vskip 0.2in

{\bf Notation and Preliminaries  }
We will use the notation $X\lesssim Y$ whenever there exists some positive constant $C$ so that $X\le C Y$. Similarly, we will use $X\sim Y$ if $X\lesssim Y \lesssim X$.
For a linear operator $A$ from Banach space $X$ to Banach space $Y$, we denote its operator norm by $\|A\|_{\mathcal{L}(X\to Y)}$.
All the constants are denoted by $C$ and they can change from line to line. We use $\varepsilon$ to denote some sufficiently small constant and it may vary from line to line.

\begin{Proposition}[Dispersive estimate of $e^{it\Delta_V}$, \cite{JAC}]\label{pro2.4}
Let $d\ge3$, ${\left\langle x \right\rangle ^\alpha }V(x)$ is a bounded operator from ${H^\eta }$ to ${H^\eta }$ for some $\alpha  > d + 4$, $\eta > 0$, with $\mathcal{F}{V}\in L^1$.
Assume also that 0 is neither an eigenvalue nor a resonance of $-\Delta_V$. Then
$${\left\| {{e^{it\Delta_V}P_c(\Delta_V)}} \right\|_{p' \to p}} \le C{\left| t \right|^{ - \frac{d}{2}\big(1 - \frac{2}{p}\big)}},
$$
where $\frac{1}{{p'}} + \frac{1}{p} = 1$, $2\le p \le \infty$.
\end{Proposition}

In our case, the spectrum of $\Delta-V$ is continuous, then
by the abstract theorem in M. Keel and T. Tao \cite{KT} and Proposition \ref{pro2.4}, one can prove:
\begin{Proposition}[Strichartz estimate]\label{1234}
For the potential $V$ in Theorem 1.1, $d\ge3$, $-\infty<t_0\le t_1<\infty$ we have
\begin{align*}
\left\| e^{it\Delta_V} f \right\|_{L_t^p L_x^q(I \times \mathbb{R}^d)} & \lesssim \|f\|_{L^2},\\
\left\|\int_{t_0}^t e^{i  (t-\tau) \Delta_V}F(\tau) \,\mathrm{d}\tau \right\|_{L^p_tL^q_x([t_0,t_1] \times \mathbb{R}^d)}¡¡& \lesssim \|F\|_{L_t^{s'}L_x^{r'}([t_0,t_1]\times\Bbb R^d)},
\end{align*}
where $(p,q)$, $(s,r)$ are Strichartz admissible with $p,q,r,s\in[2,\infty]$ namely
$$\frac{2}{p} + \frac{d}{q} = \frac{d}{2}, \mbox{  }\frac{2}{s} + \frac{d}{r} = \frac{d}{2}.$$
\end{Proposition}

\section{Resolvent estimates}\label{se2}
In this section, we will prove the resolvent estimate in Lemma \ref{as}. First, we give the basic estimate of the Green function $ G(x,y;k)$ of $(k^2-\Delta_V)^{-1}$.
\begin{Lemma}\label{q1}
If $V(x)$ satisfies (\ref{condition}), then the Green function {$ G(x,y;k)$} of $(k^2-\Delta_V)^{-1}$ satisfies the following estimate:
{for $k >0$,
 \begin{align}\label{a0}
 G(x,y;k)\le
\begin{cases}
 C(d)|x - y{|^{ - (d - 2)}},    & \mbox{  }\mbox{  }k |x - y| \le 1, \\
C(d){e^{ - k|x - y|}}|x - y{|^{ - \frac{d - 1}{2}}}|k{|^{\frac{d - 3}{2}}}, & \mbox{  }\mbox{  } k |x - y|  > 1.
 \end{cases}
\end{align}}
\end{Lemma}
\begin{proof}
It suffices to solve the following equation,
\begin{equation}\label{01}
G(x,y;k)=G_0(x,y;k)+\int G_0(x,z;k)G(z,y;k)V(z)\,\mathrm{d}z,
\end{equation}
where $G_0$ is the Green function for the resolvent of $\Delta$, and satisfies the estimate (\ref{a0}) with $G$ replaced by $G_0$.
Let
\begin{align*}
r=|x-y|, A(x,y;k)=|G(x,y;k)|r^{d-2}1_{rk \le 1}, B(x,y;k)=|G(x,y;k)|r^{\frac{d-1}{2} }k^{-\frac{d-3}{2} }1_{r k>1}e^{kr}.
\end{align*}
We prove (\ref{a0}) by using contraction mapping principle in the following Banach space:
$$X = \Big\{ G(x,y;k):{{\left\| {A(x,y;k)} \right\|}_\infty } + {{\left\| {B(x,y;k)} \right\|}_\infty } < \infty \Big\},$$
with norm defined by
$$\left\|G\right\|_X=\left\|A(x,y;k) \right\|_\infty  + \left\|B(x,y;k) \right\|_\infty.$$
Define ${\Theta}G=G_0(x,y;k)+\int G_0(x,z;k)G(z,y;k)V(z)\,\mathrm{d}z$, to prove the contraction, it suffices to verify
$$\left\| \int G_0(x,z;k)G(z,y;k)V(z)\,\mathrm{d}z\right\|_X \le \theta {\left\| {G(z,y;k)} \right\|_X},$$
for some $0<\theta<1$.
Then from (\ref{01}), we have
\begin{align*}
&\quad \int | x - y{|^{d - 2}}{1_{r k  \le 1}}\left| {{G_0}(x,z;k)G(y,z;k)V(z)} \right|\,\mathrm{d}z \\
&\le \int | x - y{|^{d - 2}}{1_{r k \le 1}}\left| {{G_0}(x,z;k)} \right|{1_{|y - z|k  \le 1}}\left| {G(y,z;k)V(z)} \right|\,\mathrm{d}z \\
&\quad + \int | x - y{|^{d - 2}}{1_{r k  \le 1}}\left| {{G_0}(x,z;k)} \right|{1_{|y - z| k  > 1}}\left| {G(y,z;k)V(z)} \right|\,\mathrm{d}z \\
&\le \int | x - y{|^{d - 2}}{1_{r|k| \le 1}}\left| {{G_0}(x,z;k)} \right|{1_{|y - z| k  \le 1}}|y - z{|^{ - (d - 2)}}A(y,z;k)|V(z)|\,\mathrm{d}z \\
& \quad + \int | x - y{|^{d - 2}}{1_{r k \le 1}}\left| {{G_0}(x,z;k)} \right||y - z{|^{ - \frac{d - 1}{2}}}|k{|^{\frac{d - 3}{2}}}{e^{ - k|y - z|}}B(y,z;k)|V(z)|\,\mathrm{d}z \\
&\lesssim\int | x - y{|^{d - 2}}{1_{r k  \le 1}}{1_{|x - z| k  \le 1}}|x - z{|^{ - (d - 2)}}|y - z{|^{ - (d - 2)}}A(y,z;k)|V(z)|\,\mathrm{d}z \\
&\quad + \int | x - y{|^{d - 2}}{1_{r k  \le 1}}{1_{|x - z| k  > 1}}{e^{ - |x - z|k}}|k{|^{\frac{d - 3}{2}}}|x - z{|^{ - \frac{d - 1}{2}}}|y - z{|^{ - (d - 2)}}A(y,z;k)|V(z)| \,\mathrm{d}z \\
&\quad + \int | x - y{|^{d - 2}}{1_{r k  \le 1}}{1_{|x - z| k  \le 1}}|x - z{|^{ - (d - 2)}}|y - z{|^{ - \frac{d - 1}{2}}}|k{|^{\frac{d - 3}{2}}}{e^{ - |y - z|k}}B(y,z;k)|V(z)|\,\mathrm{d}z \\
&\quad + \int | x - y{|^{d - 2}}{1_{r k  \le 1}}{1_{|x - z| k  > 1}}{e^{ - |x - z|k}}|k{|^{\frac{d - 3}{2}}}|x - z{|^{ - \frac{d - 1}{2}}}|y - z{|^{ - \frac{d - 1}{2}}}|k{|^{\frac{d - 3}{2}}}{e^{ - |y - z|k}}B(y,z;k)|V(z)|\,\mathrm{d}z \\
&\lesssim  \big({\left\| {A(x,y;k)} \right\|_\infty } + {\left\| {B(x,y;k)} \right\|_\infty }\big)  \int \frac{{|x - y{|^{d - 2}}}}{{|x - z{|^{d - 2}}|z - y{|^{d - 2}}}}|V(z)| \,\mathrm{d}z \\
&\lesssim  \big({\left\| {A(x,y;k)} \right\|_\infty } + {\left\| {B(x,y;k)} \right\|_\infty }\big)  \mathop {\sup }\limits_y \int {\frac{{|V(z)|}}{{|z - y{|^{d - 2}}}} \,\mathrm{d}z} .
\end{align*}
Similarly, we can prove
\begin{align*}
&\int | x - y{|^{\frac{{d - 1}}{2}}}{e^{k\left| {x - y} \right|}}{1_{r k  > 1}}{ k^{ - \frac{{d - 3}}{2}}}\left| {{G_0}(x,z;k)G(y,z;k)V(z)} \right| \,\mathrm{d}z \\
&\lesssim    ({\left\| {A(x,y;k)} \right\|_\infty } + {\left\| {B(x,y;k)} \right\|_\infty }) \mathop {\sup }\limits_y \int {\frac{{|V(z)|}}{{|z - y{|^{d - 2}}}} \,\mathrm{d}z}.
\end{align*}
Hence
$${\left\| {\int {{G_0}(x,z;k)G(y,z;k)V(z)} \,\mathrm{d}z} \right\|_X}\lesssim  {\left\| {G(y,z;k)} \right\|_X} \mathop {\sup }\limits_y \int {\frac{{|V(z)|}}{{|z - y{|^{d - 2}}}} \,\mathrm{d}z} .
$$
Therefore, by \eqref{condition}, $\Theta$ is a contraction mapping.
\end{proof}

The following estimates are immediate corollary of estimates to the Green function of resolvent and they play a prime role in our arguments.
\begin{Lemma}[Resolvent estimate]\label{as}
If $V$ satisfies (\ref{condition}), for $N>d$, we have the weighted resolvent estimates for $\lambda>0$:
\begin{align}
&{  \|(\lambda - \Delta_V)^{-1}f\|_p  \   \le C\lambda^{\frac{1}{2}(d-2-d(\frac{1}{p}+1-\frac{1}{e}))}\|f\|_e}, \label{li}\\
& \|(\lambda - \Delta_V)^{-1}\left\langle x \right\rangle^{-N}f\|_p    \le C\lambda^{-\frac{1}{p}}\|f\|_p,  \label{e}\\
& {\left\| {{{\left\langle x \right\rangle }^{ - N}}{{(\lambda  - {\Delta _V})}^{ - 1}}f} \right\|_p}    \le C{\lambda ^{ - \frac{1}{{p'}}}}{\left\| f \right\|_p}\label{kkl},
\end{align}
{where $1\le e\le p \le \infty$}.
\end{Lemma}
\begin{proof}
(\ref{kkl}) is the dual version of (\ref{e}), hence we only give the proof of (\ref{li}) and (\ref{e}).
From Young's inequality and Lemma \ref{q1}, we obtain
\begin{align*}
\|(\lambda-\Delta_V)^{-1}f\|_p\le \|G(r;\sqrt{\lambda})\|_a  {\|f\|_e} \le { \lambda^{\frac{1}{2}(d-2-d(\frac{1}{p}+1-\frac{1}{e}))}} \|f\|_e.
\end{align*}
where {$a=\frac{1}{1+ \frac{1}{p}-\frac{1}e}$}, thus we have proved (\ref{li}). Now we turn to \eqref{e}, let $c = \frac1{1+ \frac1p -\frac1a} $.
Again from Young's inequality, Lemma \ref{q1} and H\"older's inequality, we have
\begin{align*}
\|(\lambda-\Delta_V)^{-1}f\|_p&\le \|G(r;\sqrt{\lambda})\|_a\|\left\langle x \right\rangle
^{-N}  f\|_c\le C\lambda^{\frac{1}{p}}\|f\|_p\|\left\langle x \right\rangle^{-N}\|_{a'},
\end{align*}
which yields (\ref{e}).
\end{proof}

The following ``almost equivalence" Lemma is useful later.
\begin{Lemma}[Almost equivalence ]\label{34}
For any $1<p<\infty$ and $0<s<2$, $V$ in Theorem 1.1, we have
\begin{align}\label{eq2.6}
{\left\| {{{\left( { - {\Delta _V}} \right)}^{\frac{s}2}}u - {{\left( { - \Delta } \right)}^{\frac{s}2}}u} \right\|_p} \le C{\left\| u \right\|_p}.
\end{align}
\end{Lemma}
\begin{proof}
Recall the formula,
$${T^{s/2}}f = c(s){\rm{ }}T{\rm{  }}\int_0^\infty  {{\tau ^{ s/2-1}}} {\left( {\tau  + T} \right)^{ - 1}}fd\tau.
$$
Therefore one has
\begin{align*}
 {\left( { - {\Delta _V}} \right)^{s/2}}f &= c(s){\rm{ }}\left( { - {\Delta _V}} \right){\rm{  }}\int_0^\infty  {{\tau ^{ s/2-1}}} {\left( {\tau  - {\Delta _V}} \right)^{ - 1}}fd\tau  \\
 &= c(s){\rm{ }}\left( { - {\Delta _V}} \right){\rm{  }}\int_0^\infty  {{\tau ^{ s/2-1}}} \left[ {{{\left( {\tau  - {\Delta _V}} \right)}^{ - 1}} - {{\left( {\tau  - \Delta } \right)}^{ - 1}}} \right]fd\tau  + c(s){\rm{ }}\left( { - {\Delta _V}} \right){\rm{  }}\int_0^\infty  {{\tau ^{s/2-1}}} {\left( {\tau  - \Delta } \right)^{ - 1}}fd\tau  \\
 &= c(s){\rm{ }}{\Delta _V}{\rm{  }}\int_0^\infty  {{\tau ^{ s/2-1}}} {\left( {\tau  - {\Delta _V}} \right)^{ - 1}}V{\left( {\tau  - \Delta } \right)^{ - 1}}fd\tau  + c(s){\rm{ }}\left( { - {\Delta _V}} \right){\rm{  }}\int_0^\infty  {{\tau ^{s/2-1}}} {\left( {\tau  - \Delta } \right)^{ - 1}}fd\tau  \\
&= - c(s)\int_0^\infty  {{\tau ^{s/2-1}}} \left( {\tau  - {\Delta _V}} \right){\left( {\tau  - {\Delta _V}} \right)^{ - 1}}V{\left( {\tau  - \Delta } \right)^{ - 1}}fd\tau  \\
&+ c(s)\int_0^\infty  {{\tau ^{s/2}}} {\left( {\tau  - {\Delta _V}} \right)^{ - 1}}V{\left( {\tau  - \Delta } \right)^{ - 1}}fd\tau
+ c(s){\rm{ }}\left( { - \Delta } \right)\int_0^\infty  {{\tau ^{s/2-1}}} {\left( {\tau  - \Delta } \right)^{ - 1}}fd\tau\\
&+ c(s)\int_0^\infty  {{\tau ^{s/2-1}}} V{\left( {\tau  - \Delta } \right)^{ - 1}}fd\tau  \\
&= {\left( { - \Delta } \right)^{s/2}}f + c(s)\int_0^\infty  {{\tau ^{s/2}}} {\left( {\tau  - {\Delta _V}} \right)^{ - 1}}V{\left( {\tau  - \Delta } \right)^{ - 1}}fd\tau.
\end{align*}
Therefore we have proved
$$(-\Delta_V)^{\frac{s}2}=(-\Delta)^{\frac{s}2}+\int^{\infty}_0\lambda^{s/2}(\lambda-\Delta_V)^{-1}V(\lambda-\Delta)^{-1}\,\mathrm{d}\lambda.$$
Then from Lemma \ref{as}, we have
$$\left\|\int^{\infty}_1\lambda^{s/2}(\lambda-\Delta_V)^{-1}V(\lambda-\Delta)^{-1}u \,\mathrm{d}\lambda\right\|_p
\le C \|u\|_p \|V\|_{\infty}\int^{\infty}_1\lambda^{s/2-2}\,\mathrm{d}\lambda,
$$
and
\begin{align*}
& {\left\| {\int_0^1 {{\lambda ^{s/2}}} {{(\lambda  - {\Delta _V})}^{ - 1}}V{{(\lambda  - \Delta )}^{ - 1}}u\,\mathrm{d}\lambda } \right\|_p}\\
\le  & \  {\left\| u \right\|_p} {\left\| {{{\left\langle x \right\rangle }^{2N}}V} \right\|_\infty }\int_0^1 {{\lambda ^{s/2-1}}} {\left\| {{{(\lambda  - {\Delta _V})}^{ - 1}}{{\left\langle x \right\rangle }^{ - N}}} \right\|_{p \to p}}{\left\| {{{\left\langle x \right\rangle }^{ - N}}{{(\lambda  - \Delta )}^{ - 1}}} \right\|_{p \to p}}\mathrm{d}\lambda \\
\le & \  \|u\|_p \int^{1}_{0}\lambda^{\frac{s}2-1}\,\mathrm{d}\lambda\\
\le & \  C\|u\|_p,
\end{align*}
thus \eqref{eq2.6} follows.
\end{proof}

\section{The commutator operator $|J_V(t)|^s$ and proof of Theorem \ref{th5}}\label{se3}
This section is devoted to the proof of Theorem \ref{th5}. In subsection \ref{subse3.1}, we will estimate the commutator $|J_V(t)|^s u$. In subsection \ref{subse3.2}, we give the decay estimate, and thus scattering follows easily from the decay estimate.
\subsection{The commutator operator $|J_V(t)|^s$}\label{subse3.1}
In \cite{SVN}, the authors proposed the following commutator operator: $\forall\, 0 < s <2$,
$$|J_V(t)|^s u=M(t)(-t^2\Delta_V)^{\frac{s}2}M(-t)u,$$
where $M(t)=e^{i\frac{|x|^2}{4t}}$.
Moreover, they proved that $|J_V|^su$ satisfies
$$i\partial_t|J_V|^su +  \Delta_V |J_V|^su  - it^{s-1}M(t)A(s)M(-t)u  + \lambda |J_V|^s\big(|u|^{p-1} u\big)=0,
$$
where
\begin{equation}\label{r}
A(s)=c(s)\int^{\infty}_{0}\tau^{\frac{s}2}(\tau-\Delta_V)^{-1}(2V+x \cdot {\nabla_x}V)(\tau-\Delta_V)^{-1}\,\mathrm{d}\tau.
\end{equation}
Let $s=\big(\frac{d}2\big( 1 -\frac{1}{p}\big)\big)^{+}$.
By Proposition \ref{1234}, we have
\begin{equation} \label{fd}
\big\||J_V|^su\big\|_{L^{\infty}_tL^2_x}\le C\big\||J_V(h)|^su(h)\big\|_2+C\|t^{s-1}A(s)M(-t)u\|_{L^r_tL^q_x}+C\big\||J_V|^s\big( |u|^{p-1}u\big)\big\|_{L^m_tL^k_x},
\end{equation}
where $(r',q')$ and $(m',k')$ are admissible pair.
Choose $1<m,r, q,k<2$ such that
\begin{equation}\label{sd}
\frac{s}{2}-\frac{1}{q}+\frac{1}{2p}>0,
\end{equation}
\begin{equation}\label{wr}
 p>1+\frac{2m-1}{\frac{1}{2}dm},
\end{equation}
which are possible since $d>2$, $p>1+\frac{2}{d}$.

Therefore, in order to apply continuity method, we have to bound $\|t^{s-1}A(s)M(-t)u\|_{L^r_tL^q_x}$ and $\big\||J_V|^s|u|^{p-1}u\big\|_{L^m_tL^k_x}$ by $\big\||J_V|^su\big\|_{L^{\infty}_tL^2_x}$.
Preliminarily, we present some properties of the commutator operator, especially ``Sobolev embedding theorem".

\begin{Lemma}\label{qw}
For $u\in H^s$, $s=\big(\frac{d}2\big( 1 -\frac{1}{p}\big)\big)^{+}$, there exists some $0<\eta<1$ such that
$$\|u\|_{2p}\le C\big\|(-\Delta_V)^{\frac{s}2}u\big\|_2+C\big\|(-\Delta_V)^{\frac{s}2}\big\|^{\eta}_2\|u\|^{1-\eta}_2,
$$
\end{Lemma}
\begin{proof}
The proof is divided into four steps.

{\bf{Step 1.}}
We reduce the problem to $u\in C_c^{\infty}$ by density arguments. We claim that for $u\in H^s$, there exists a sequence of functions $u_n\in C^{\infty}_c$ such that $\big\|(-\Delta_V)^{\frac{s}2}(u-u_n)\big\|_2\to 0.$
In fact, since $u\in H^s$, there exists $u_n\in C^{\infty}_c$ such that $\big\|(-\Delta)^{\frac{s}2}(u_n-u)\big\|_2+\|u_n-u\|_2\to 0$.
Lemma \ref{34} implies
$$\big\|(-\Delta_V)^{\frac{s}2}(u-u_n)-(-\Delta)^{\frac{s}2}(u-u_n)\big\|_2\le C\|u_n-u\|_2,$$
by which our claim follows.
Therefore, we can assume $u\in C^{\infty}_c$ without loss of generality.

{\bf{Step 2}}. Reduction to the estimate of $(-\Delta_V)^{-\frac{s}2}$.
Formally, we have
$$(-\Delta_V)^{-\frac{s}2}=c(s)\int^{\infty}_{0}{\lambda}^{-\frac{s}2}(\lambda-\Delta_V)^{-1}\,\mathrm{d}\lambda.$$
Precisely, it converges strongly in $L^{2p}$ for $u\in C^{\infty}_{c}$. In fact, \eqref{li} implies
$$\int^{\infty}_1\|{\lambda}^{-\frac{s}2}(\lambda-\Delta_V)^{-1}u\|_{2p}\,\mathrm{d}\lambda
\le \int^{\infty}_{1}\lambda^{-\frac{s}2-1}\|u\|_{2p}\,\mathrm{d}\lambda,
$$
and
$$\int^{1}_0\|{\lambda}^{-\frac{s}2}(\lambda-\Delta_V)^{-1}u\|_{2p}\,\mathrm{d}\lambda
\le \lambda\|u\|_{\alpha} \int^{1}_{0}{\lambda}^{-\frac{s}2+\frac{1}{2}(d-2-d(\frac{1}{2p}+1-\frac{1}{\alpha}))}\,\mathrm{d}\lambda
\le C\|u\|_{\alpha},
$$ for $\frac{1}{\alpha}>\frac{1}{2p}+\frac{s}{d}$.
Therefore it suffices to prove
$$\big\|(-\Delta_V)^{-\frac{s}2}u\big\|_{2p}\le C\|u\|_2+C\|u\|^{\eta}_2\big\|(-\Delta_V)^{-\frac{s}2}u\big\|^{1-\eta}_2.
$$

{\bf{Step 3}}. High energy estimate. From (\ref{li}), we obtain

$$\left\|\int^{\infty}_1{\lambda}^{-\frac{s}2}(\lambda-\Delta_V)^{-1}u\,\mathrm{d}\lambda\right\|_{2p}
\le \|u\|_2 \int^{\infty}_{1}{\lambda}^{-\frac{s}2+\frac{1}{2}(d-2-d(\frac{1}{2p}+\frac{1}{2}))}\,\mathrm{d}\lambda\le C\|u\|_2,
$$
where we have used $s=\big(\frac{d}2\big(1-\frac{1}{p}\big)\big)^+$.

{\bf{Step 4}}. Low energy estimate.
Define $\frac{1}{2}-\frac{1}{\gamma}=\frac{s}{d}$, $\mu>\gamma$, $\frac{1- \eta}{\mu}+\frac{\eta}{2}=\frac{1}{2p}$.
Then H\"older's inequality indicates
$$\left\|\int^{1}_0\lambda^{-\frac{s}2}(\lambda-\Delta_V)^{-1}u\,\mathrm{d}\lambda
\right\|_{2p} \le \left\|\int^1_0\lambda^{-\frac{s}2}(\lambda-\Delta_V)^{-1}u\,\mathrm{d}\lambda\right\|^{1- \eta}_{\mu}\left\|\int^1_0
\lambda^{-\frac{s}2}(\lambda-\Delta_V)^{-1}u\,\mathrm{d}\lambda\right\|^{\eta}_{2}.
$$
\eqref{li} gives
$$
\left\|\int^1_0\lambda^{-\frac{s}2}(\lambda-\Delta_V)^{-1}u\,\mathrm{d}\lambda\right\|_{\mu}
\le \|u\|_2 \int^1_0\lambda^{-\frac{s}2+\frac{1}{2}(d-2-d(\frac{1}{\mu}+\frac{1}{2}))}\,\mathrm{d}\lambda\le C\|u\|_2,
$$
where we have used $\mu>\gamma$.
Again by \eqref{li}, we obtain
\begin{align*}
\left\|\int^1_0 \lambda^{-\frac{s}2}(\lambda-\Delta_V)^{-1}u\,\mathrm{d}\lambda\right\|_{2}
&\le \left\|\int^{\infty}_0\lambda^{-\frac{s}2}(\lambda-\Delta_V)^{-1}u\,\mathrm{d}\lambda\right\|_{2}
+\|u\|_2 \int^{\infty}_1\lambda^{-\frac{s}2-1}\,\mathrm{d}\lambda\\
&\le \big\|(-\Delta_V)^{-\frac{s}2}u\big\|_2+C\|u\|_2.
\end{align*}
Assembling all the estimates above, we finish the proof.
\end{proof}

As a direct consequence of Lemma \ref{qw}, we easily obtain the following lemma.
\begin{Lemma}\label{15}
Taking $s=\big(\frac{d}2\big( 1 -\frac{1}{p}\big)\big)^+$, $\mu=(2p)^+$, then for $s_0=\big(\frac{d}2\big( 1 -\frac{1}{ p}\big)\big)^{-}$, it holds
$$\|u\|_{2p}\le Ct^{-s_0}\Big(\big\||J_V|^{\frac{s}2}u\big\|_{L^{\infty}_tL^2_x}+\big\||J_V|^{\frac{s}2}u\big\|^{\eta}_{L^{\infty}_tL^2_x}\|u\|^{1-\eta}_2\Big).
$$
\end{Lemma}

\begin{Lemma}\label{16}
$$\|A(s)u\|_q\le C\|u\|_{2p}.$$
\end{Lemma}
\begin{proof}
Let $V_1 = 2V+x\cdot\nabla_x V$ and $ w = \frac{2pq}{2p-q}$, Lemma \ref{as} and (\ref{r}) give
\begin{align*}
\|A(s)u\|_q
& \ \le c(s)\left\|\int^{\infty}_{1}\tau^{\frac{s}2}(\tau-\Delta_V)^{-1} V_1(\tau-\Delta_V)^{-1}u\,\mathrm{d}\tau\right\|_q+c(s)\left\|\int^{1}_{0}\tau^{\frac{s}2}
(\tau-\Delta_V)^{-1} V_1(\tau-\Delta_V)^{-1}u\,\mathrm{d}\tau\right\|_q\\
& \ \le \|u\|_{2p} \  c(s)\int^{\infty}_{1}\tau^{\frac{s}2}\|(\tau-\Delta_V)^{-1}\|_{p\to p}\|V_1\|_w\|(\tau-\Delta_V)^{-1}\|_{2p\to 2p}\,\mathrm{d}\tau\\
& \ \quad  + \|u\|_{2p} \
c(s)\int^{1}_{0}\tau^{\frac{s}2}\|(\tau-\Delta_V)^{-1}\left\langle x \right\rangle^{-N}\|_{p\to p}\|\left\langle x \right\rangle^{2N}V_1\|_w\|\left\langle x \right\rangle^{-N}(\tau-\Delta_V)^{-1}\|_{2p\to2p}\,\mathrm{d}\tau\\
&\ \le  C \|u\|_{2p} \int^{\infty}_{1}\tau^{\frac{s}2-2}\,\mathrm{d}\tau +C \|u\|_{2p} \int^{1}_{0}\tau^{\frac{s}2-\frac{1}{p}+\frac{1}{2p}-1}\,\mathrm{d}\tau \\
&\ \le  C\|u\|_{2p}.
\end{align*}
\end{proof}

Lemma \ref{16} combined with Lemma \ref{15} implies that if $r$ is appropriately chosen, then
\begin{Corollary}\label{df}
For any $u: \mathbb{R} \times  \mathbb{R}^d \to \mathbb{C}$, we have
\begin{equation*}
\|t^{s-1}M(t)A(s)M(-t)u\|_{L^r_tL^q_x}\le C(h)\Big(\big\||J_V|^su\big\|_{L^{\infty}_tL^2_x}+\|u\|_{L_t^\infty L_x^2}\Big).
\end{equation*}
\end{Corollary}
Now we turn to the nonlinear term.
\begin{Lemma}\label{fg}
There exists $0<\theta<1$ such that
\begin{align}\label{eq3.5}
 \big\||J_V|^s(|u|^{p-1}u)\big\|_{L^m_tL^k_x}
&\le C\Big(\big\||J_V|^su\big\|_{L^{\infty}_tL^2_{x}}+\|u\|_{L_t^\infty L_x^2}\Big)\Big(\big\||J_V|^su\big\|^{\theta}_{L^{\infty}_tL^2_x}
+\big\||J_V|^su\big\|^{\theta\eta}_{L^{\infty}_tL^2_x}\|u\|^{(1-\eta)\theta}_{L_t^\infty L_x^2}\Big)^{p-1}\notag\\
&\quad +C(h)\Big(\big\||J_V|^su\big\|^{\theta}_{L^{\infty}_tL^2_x}+\big\||J_V|^su\big\|^{\theta\eta}_{L^{\infty}_tL^2_x}\|u\|^{(1-\eta)\theta}_{L_t^\infty L_x^2}\Big)^p
\end{align}
\end{Lemma}
\begin{proof}
Let $\tilde{m} = \frac{2k}{2-k}$, by Lemma 3.4 in \cite{JTG} and Lemma \ref{34}, we have
\begin{align*}
\big\||J_V|^s(|u|^{p-1}u)\big\|_{L^m_tL^k_x}
&\le\big\|(-\Delta)^{\frac{s}2}M(-t)(|u|^{p-1}u)\big\|_{L^m_tL^k_x}+\big\||u|^{p-1}u\big\|_{L^m_tL^k_x}\\
&
\le\big\|(-\Delta)^{\frac{s}2}\big(|M(-t)u|^{p-1}M(-t)u\big)\big\|_{L^m_tL^k_x}+\big\|\|u\|^p_{L^{kp}_x}\big\|_{L_t^m}\\
&
\le\big\|(-\Delta)^{\frac{s}2}(M(-t)u)\big\|_{L^{\infty}_tL^2_x}\big\||u|^{p-1}\big\|_{L^m_tL^{\tilde m}_x}+\big\|\|u\|^p_{L^{kp}_x}\big\|_{L_t^m}\\
&
\le C\Big(\big\||J_V|^su\big\|_{L^{\infty}_tL^2_x}+\|u\|_2\Big)\big\||u|^{p-1}\big\|_{L^m_tL^{\tilde m}_x}+\big\|\|u\|^p_{L^{kp}_x}\big\|_{L_t^m}.
\end{align*}
First, we consider $\big\|\|u\|^p_{L^{kp}_x}\big\|_{L_t^m}$. H\"older's inequality and Lemma \ref{15} show
\begin{align*}
\|u\|_{L^{kp}_x}&\le\|u\|^{\theta}_{L^{2p}_x}\|u\|^{1-\theta}_2\\
&\le t^{-d(\frac{1}{2}-\frac{1}{kp})+\varepsilon}\Big(\big\||J_V|^su\big\|^{\theta}_{L^{\infty}_tL^2_x}+\big\||J_V|^su\big\|^
{\theta\eta}_{L^{\infty}_tL^2_x}\|u\|^{(1-\eta)\theta}_2\Big) \|u\|^{1-\theta}_2,
\end{align*}
where $\varepsilon$ is sufficiently small.
Since $(m',k')$ is an admissile pair, we have $k=\frac{1}{\frac{1}{2}+\frac{2}{d}(1-\frac{1}{m})},$ then
\begin{align*}
\big\|\|u\|^p_{L^{kp}_x}\big\|_{L^m}
&\le
\left(\int^{\infty}_ht^{\varepsilon-dmp\Big(\frac{1}{2}-\frac{\frac{1}{2}+\frac{2}{d}(1-\frac{1}{m})}{p}\Big)}\,\mathrm{d}t \right)^{\frac{1}{m}}\Big(\big\||J_V|^su\big\|
^{\theta}_{L^{\infty}_tL^2_x}+\big\||J_V|^su\big\|^{\theta\eta}_{L^{\infty}_tL^2_x}\|u\|^{(1-\eta)\theta}_2\Big)^p\\
&\le C(h)\Big(\big\||J_V|^su\big\|^{\theta}_{L^{\infty}_tL^2_x}+\big\||J_V|^su\big\|^{\theta\eta}_{L^{\infty}_tL^2_x}\|u\|^{(1-\eta)\theta}_2\Big)^p,
\end{align*}
where we have used (\ref{wr}).
Second, we estimate $\big\||u|^{p-1}\big\|_{L^m_tL^{\tilde{m} }_x}$.
Similar arguments as the above estimates give
\begin{align*}
\big\||u|^{p-1}\big\|_{L^m_tL^{\tilde{m}}_x}
&\le \left(\int^{\infty}_h t^{\varepsilon-d(\frac{1}{2}-\frac{\frac{1}{2}-\frac{1}{d}\left(\frac{n}{2}-\frac{2}{m'}\right)}{p-1})(p-1)m} \,\mathrm{d}t \right)^{\frac{1}{m}}
\Big(\big\||J_V|^su\big\|^{\theta}_{L^{\infty}_tL^2_x}+\big\||J_V|^su\big\|^{\theta\eta}_{L^{\infty}_tL^2_x}\|u\|^{(1-\eta)\theta}_2\Big)^{p-1}\\
&\le C(h)\Big(\big\||J_V|^su\big\|^{\theta}_{L^{\infty}_tL^2_x}+\big\||J_V|^su\big\|^{\theta\eta}_{L^{\infty}_tL^2_x}\|u\|^{(1-\eta)\theta}_2\Big)^{p-1},
\end{align*}
where again we have used (\ref{wr}).
Combining the estimates together, we obtain \eqref{eq3.5}.
\end{proof}

\subsection{Proof of Theorem \ref{th5}}\label{subse3.2}
Now we are ready to prove Theorem \ref{th5}. By (\ref{fd}), Corollary \ref{df} and Lemma \ref{fg}, we have
\begin{align*}
\big\||J_V|^su\big\|_{L^{\infty}_tL^2_x}  & \le    C\Big(\big\||J_V|^su\big\|_{L^{\infty}_tL^2_{x}}+\|u\|_2\Big)\Big(\big\||J_V|^su\big\|^{\theta}_{L^{\infty}_tL^2_x}
+\big\||J_V|^su\big\|^{\theta\eta}_{L^{\infty}_tL^2_x}\|u\|^{(1-\eta)\theta}_2\Big)^{p-1}\\
& \quad  +C(h)\Big(\big\||J_V|^su\big\|^{\theta}_{L^{\infty}_tL^2_x}+\big\||J_V|^su\big\|^{\theta\eta}_{L^{\infty}_tL^2_x}\|u\|^{(1-\eta)\theta}_2\Big)^p\\
&\quad  +C\|u_0\|_{\Sigma}+C(h)\big\||J_V|^su\big\|_{L^{\infty}_tL^2_{x}}.
\end{align*}
Since $\mathop {\lim }\limits_{h \to \infty }C(h)=0$, by standard continuity argument, we obtain for sufficiently large $h$ and sufficiently small $\|u_0\|_{\Sigma}$
$$
\big\||J_V|^su\big\|_{L^{\infty}_tL^2_x}\le C.
$$
We then have the decay estimate \eqref{eq1.3} by Lemma 3.2.

Finally, we give the proof of scattering as a consequence of \eqref{eq1.3}.
From Duhamel's principle, it suffices to prove
$$\left\|\int^{\infty}_he^{-is\Delta_V}\left(\lambda|u|^{p-1}u\right)(s)\,\mathrm{d}s\right\|_{H^1}\le C.
$$
Lemma \ref{34} and Strichartz estimate yield
\begin{align*}
\left\|\int^{\infty}_he^{-is\Delta_V}\left(\lambda |u|^{p-1}u\right)(s)\,\mathrm{d}s\right\|_{H^1}
&\lesssim \left\|\int^{\infty}_he^{-is\Delta_V}\left(\lambda |u|^{p-1}u\right)(s)\,\mathrm{d}s\right\|_2\\
&+\left\|\int^{\infty}_h(-\Delta_V)^{\frac12}e^{-is\Delta_V}\left(\lambda |u|^{p-1}u\right)(s)\,\mathrm{d}s\right\|_2\\
&
\lesssim \big\|\lambda |u|^{p-1}u\big\|_{L^m_tL^k_x}+\big\|(-\Delta_V)^{\frac12}\left(\lambda |u|^{p-1}u\right)\big\|_{L^m_tL^k_x}\\
&
\lesssim \big\|\lambda |u|^{p-1}u\big\|_{L^m_tL^k_x}+\big\|(-\Delta)^{\frac12}\left(\lambda |u|^{p-1}u\right)\big\|_{L^m_tL^k_x}\\
&
\lesssim \big\||u|^{p-1}u\big\|_{L^m_tL^k_x}+\|u\|_{L^{\infty}_tH^1_x}\big\||u|^{p-1}\big\|_{L^m_tL^{\tilde{m}}_x},
\end{align*}
where $\tilde{m} = \frac{2k}{2-k}$.
The argument of the proof of Lemma \ref{fg} implies
$$\big\||u|^{p-1}u\big\|_{L^m_tL^k_x}+C\big\||u|^{p-1}\big\|_{L^m_tL^{\bar m}_x}\le C.
$$
We then define
\begin{equation*}
u_+ = e^{-ih\Delta_V} u_0 - i \int_h^\infty e^{-i\tau \Delta_V}\big(\lambda|u|^{p-1}u\big)(\tau)\,\mathrm{d}\tau,
\end{equation*}
and
\begin{equation*}
\|e^{-it\Delta_V } u(t) - u_+\|_{H^1} \to 0, \text{ as } t \to \infty,
\end{equation*}
 which yields the scattering.

\begin{tabular}{@{}r@{}p{16cm}@{}}
&Ze Li, {\small{Wu Wen-Tsun Key Laboratory of Mathematics, Chinese Academy of Sciences and Department of Mathematics, University of Science and Technology of China, Hefei 230026, \ Anhui, \ China}}.\\
&E-mail: lize@mail.ustc.edu.cn;\\
&Lifeng Zhao, {\small{Wu Wen-Tsun Key Laboratory of Mathematics, Chinese Academy of Sciences
 and Department of Mathematics, University of Science and Technology of China, Hefei 230026, \ Anhui, \ China}}.\\
&E-mail: zhaolf@ustc.edu.cn.
\end{tabular}

\end{document}